\documentclass[smallextended]{svjour3}
\usepackage{latexsym, amsmath, amsfonts, amscd, amssymb, verbatim}
\usepackage{xcolor}
\newcommand{\co}{\operatornamewithlimits{co}}

\newcommand{\argmax}{\operatornamewithlimits{argmax}}

\newtheorem{Remark}{Remark}

\begin{document}
\title{Algorithms for finding global and local equilibrium points of Nash-Cournot equilibrium models involving concave cost}
\author{Le Dung Muu \and Nguyen Van Quy}
\institute{Le Dung Muu, Corresponding author \at
            TIMAS, Thang Long University \\
            ldmuu@math.ac.vn
			\and
			Nguyen Van Quy \at
            Department of Mathematics, Institute of Finance, Dong Ngac, Tu Liem, Hanoi, Vietnam\\
            quynv2002@yahoo.com
}
\date{Received: date / Accepted: date}


\maketitle

\begin{abstract} We consider Nash-Cournot oligopolistic equilibrium models involving separable concave cost functions. In contrast to the models with linear and
convex cost functions, in these models a local equilibrium point may not be a global one. We propose algorithms for finding global and local equilibrium points
 for the models having separable concave cost functions. The proposed algorithms use the convex envelope of a separable concave cost function over boxes
to approximate a concave cost model with an affine cost one. The latter is equivalent to a strongly convex quadratic program that can be solved efficiently.  To obtain better approximate solutions the algorithms use an adaptive rectangular bisection which is performed only in the space of concave variables Computational results on a lot number of randomly generated data show that the proposed algorithm for global equilibrium point  are efficient for the models
with moderate number of concave cost functions while the algorithm for local equilibrium point can solve efficiently the models with much larger size.
 \subclass{$49$M$37$ \and $90$C$26$ \and $65$K$15$}
\end{abstract}

\section{Introduction} The Nash- Cournot oligopolistic market  model is one of fundamental models in economics that has been earned attention of many authors, see e.g. \cite{Au1,Fa1,Fa2,Fu2,Ko1,Ku1,Mur1} and the references cited therein. In this model it is assumed that there are $N$-firms producing a common homogeneous commodity. Each firm $i$  has a  strategy set $D_i \subset\Bbb R_{+}$
  and a profit function $f_i$ defined on the strategy set $D:= D_1 \times\cdots\times D_N$ of the model. Let $x_i\in D_i$ be a corresponding production level
  of firm $i$. Actually,  each firm seeks to  maximize its profit by choosing the corresponding production level under the presumption that  the production of the
   other firms are parametric input.  A commonly used approach to this model is based upon the famous Nash equilibrium concept.\\
\indent We recall that a point (strategy) $x^* = (x^*_1,\ldots,x^*_N) \in D$ is said to be a Nash equilibrium point of this Nash-Cournot
     oligopolistic market  model if  $$ f_i(x^*) \geq f_i(x^*[x_i]) \ \forall x_i \in D_i, \ \forall i,$$
    where the vector $x^*[x_i]$ is obtained from $x^*$ by replacing $x^*_i$ with $x_i$.\\
\indent In the  linear Nash-Cournot model the profit function of firm $i$ is given by
\begin{equation}\label{1}
f_i(x) = (\alpha - \beta\sum_{j=1}^{N} x_j)x_i - h_i(x_i) \ (i=1,\ldots,N),
\end{equation}
 where $ \beta > 0$, $ \alpha > 0$   and, for every $i$,  the cost function $h_i$ is affine that depends only on the quantity $x_i$ of firm $i$. In this linear case, it has been shown that (see e.g. \cite{Ko1}) the model has a unique Nash equilibrium  point which is the unique solution of a strongly convex quadratic program. In the case $h_i$ is differentiable convex, the problem of finding a Nash equilibrium point can be formulated as a monotone variational inequality \cite{Fa1,Mur1} which can be solved by available methods for the monotone variational inequality.\\
\indent In some practical applications,  the cost for production of a unit commodity decreases as the quantity of the production gets larger. The cost function then is concave rather than convex. Nash-Cournot oligopolistic  models with concave cost functions are considered in recent paper by Bigi and Passacantando  in \cite{Bi2}. For these models, as it is shown \cite{MHQ} that the problem can be formulated as a mixed variational inequality
  of the form
$$\text{Find}\ x^* \in D:  \langle F(x^*), x-x^*\rangle + \varphi (x) - \varphi(x^*) \geq 0\ \forall x\in D.$$
In this problem $F$ is not monotone and $\varphi$ may not be convex, and therefore the existing methods for the monotone variational inequality cannot be
applied. In  \cite{MHQ} an algorithm is proposed for  finding a global equilibrium point of the model when some of the cost
functions are piecewise linear concave. However the algorithm there is efficient only when the number of the piecewise linear concave
cost functions is relatively small. In \cite{QM1} a proximal point method was described for finding a stationary point of the model.
However a stationary point may not be a global, even not a local equilibrium point.\\
\indent In this paper we continue our work in \cite{MHQ} and \cite{QM1} by  considering Nash-Cournot models, where some of the cost functions are separable
concave, the remaining costs are affine. Namely we approximate the model with concave cost functions by   piecewise linear concave cost models that can be solved by an existing Search-and-Check algorithm in \cite{MHQ}.  Thanks to the fact that the strategy set is a rectangle (box) and the cost functions are separable increasing, the model has particular features that can be employed to develop efficient algorithms for solving it.  We propose two algorithms: the
first  one    is a search-check-branch procedure that approximates the model with concave cost functions by the models with piecewise linear concave functions. Thanks to the affine property of the price function and separability of the concave cost function the latter models can be equivalently formulated as a strongly convex quadratic problem. In order to obtain better approximate solutions the algorithm use an adaptive rectangular bisection which is performed
only in the space of the concave variables. The computational results on a lot number of randomly generated data show that this algorithm are efficient for models with a medium number ($\leq 40)$ of the firms having concave cost functions, the number of total variables may be much larger. In order to solve the models with larger number of the firms having concave cost functions we use again the convex envelope of a concave function over a box  to
develop an algorithm for obtaining a local equilibrium point.\\
\indent The remaining  part of the paper is organized as follows. In the next section we define a gap function that can serve as a stoping criterion for the algorithms. The third section is devoted to description of the algorithms and analysis of their convergence. We close the paper with some computational
results and experiences.
\section {A Gap Function as a Stoping Criterion}
In this section, we define a gap function for  Nash-Cournot models involving concave cost functions. This gap function will serve as a stoping criterion for checking whether a point is equilibrium or not. To be precise,  we consider the Nash- Cournot oligopolistic market model presented above under the assumption that  each profit function $f_j$ is defined by \eqref{1} where $h_j$, $j=1,\ldots,n$ with $n\leq N$ is increasing concave while $h_i$ with
$i>n$ is increasing affine. This assumption is motivated by the fact that for some firms the cost  consists of both the production and transportation costs, while for the  other  ones, the production need not to transport. In practice the transportation cost function is concave (see the example in \cite{Bi2}).\\
\indent First, we define the bifunction $\phi$ by taking
\begin{equation}\label{4}
 \phi(x, y) := \langle \tilde{B_1}x - a, y - x \rangle
 + y^TB_1y - x^TB_1x +   h(y) -
 h(x)\end{equation}
where $$a :=(\alpha, \alpha,\ldots,\alpha)^T,$$
$$B_1 := \begin{pmatrix} \beta &0&0&\ldots&0\cr
                  0&\beta &0& \ldots&0 \cr
                  \ldots&\ldots&\ldots&\ldots&\ldots\cr
                  0&0&0&0&\beta \end{pmatrix}, \
\tilde{B_1} := \begin{pmatrix} 0&\beta &\beta &\ldots&\beta \cr
                                        \beta &0&\beta &\ldots&\beta \cr
                                         \ldots&\ldots&\ldots&\ldots&\ldots\cr
                                         \beta &\beta&\beta&\ldots&0\end{pmatrix},$$
and we suppose that
$$h(x) := \sum_{i = 1}^N h_i(x_i).$$
Then the problem of finding an equilibrium point for the model can be formulated as a  mixed variational inequality problem MV$(D)$ of the form (see e.g.\cite {MHQ})
$$ \begin{cases} \text{find a point} \  x\in D \ \text{such that}\\
\Phi(x, y):= \langle \tilde{B_1}x   - a, y - x \rangle + \varphi(y)
- \varphi(x)  \geq 0 \ \forall y\in D,
\end{cases}\eqno MV(D)$$
where $\varphi (y):=  y^TB_1y  + h(y)$, $\varphi(x):=  x^TB_1x  + h(x)$. Clearly, $\varphi$  is a DC separable function if each $h_i$ is concave, in particular case, if each $h_i$ is affine, then $\varphi$ is a separable strongly quadratic convex function. In the latter case every local equilibrium point is global one and we have the following lemma.
\begin{lemma}\cite{Ko1,MHQ}\label{2.1a}
Suppose that the cost function $h$ is  affine(classical model) given as
$h(x) := \mu^Tx + \xi$. Then variational inequality $MV(D)$ can be
equivalently formulated as the convex quadratic programming problem
$$\min\{ x^T (2B_1 +\tilde{B_1}) x + (\mu -a)^T: x\in D\}.$$
\end{lemma}
\indent Gap functions are commonly used to determine stoping rules in optimization, variational inequality and equilibrium problems as well as to reformulate them  as a mathematical programming  problem. Following this idea,  we now  define a  gap function for the   Nash-Cournot equilibrium  models with  separable concave cost functions.  Namely, for Problem MV$(D)$ we define a gap function by taking, for each $x\in D$,
\begin{equation}\label{gap1}
g(x):=-\min\{\Phi(x,y): y\in D\}.
\end{equation}
\begin{lemma} \label{3.1} Suppose that cost function $h_i$ is continuous on $D_i$ for all $i=1, 2, \ldots, N$.  Then\\
\indent (i) The function $g(x)$ is well defined, continuous  and $g(x) \geq 0 \ \forall x\in D$;\\
\indent (ii) A point $x^*\in D$ is  equilibrium for the model if only if $g(x^*)=0$.
\end{lemma}
{\it Proof.}
This lemma can be derived from Theorem 2.1 in \cite{Fu2}. Here we give a direct proof for MV$(D)$.\\
\indent (i) Since $D$ is compact and, for each $x\in D$, $\Phi(x, .)$ is continuous on $D$, $\Phi(x,.)$ attains its minimum
 on $D$. Further, from property $\Phi(x,x)=0$, it follows that $g(x) \geq 0$ for every $x\in D$.\\
\indent (ii) Suppose that $x^* \in D$ is an equilibrium point, then
$$\Phi(x^*,y) \geq 0 \ \text{for all} \ y\in D,$$
which  implies $g(x^*) \leq 0.$ Hence  $g(x^*)= 0$. Conversely, if
$g(x^*)=0$, then from the definition of $g(x^*)$  one has
$\Phi(x^*,y)\geq 0$ for all $y\in D$, that means that  $x^*$ is a
equilibrium point of the model. $\hfill \square$\\
\indent Motivated by this lemma, we call a point $x_\epsilon$ an $\epsilon$-equilibrium point if $g(x_\epsilon) \leq \epsilon$.\\
 We rewrite the bifunction $\Phi$ as
$$\Phi(x, y) =  \langle \tilde{B_1}x - a, y - x \rangle
 +  \beta \sum_{i=1}^N  y^2_i +  \sum_{i=1}^N h_i(y_i) -
  \beta\sum_{i=1}^N
  x^2_i   - \sum_{i=1}^N h_i(x_i),$$
the gap function $g$ then can be rewritten as
\begin{equation}\label{g2}
g(x)= - \min_{y\in D}\Big \{ \langle \tilde{B_1}x - \alpha, y - x
\rangle+  \beta\sum_{i=1}^N  y^2_i +\sum_{i=1}^N h_i(y_i)\Big \} +
\beta \sum_{i=1}^Nx_i^2 +\sum_{i=1}^N h_i(x_i)
\end{equation}
Since $D$ is the box of the form
$$D:= \{ x^T = (x_1,\ldots,x_N): 0\leq l_i \leq x_i \leq u_i, \ i=1,\ldots,N\}$$
we can further write $g(x)$ as
\begin{equation}\label{g3}
\begin{aligned}
g(x)= - &\sum_{i=1}^N \min_{l_i\leq y_i\leq u_i}\Big \{  (
\tilde{B_1}x
- \alpha)_i (y_i-x_i) + \beta  y^2_i + h_i(y_i)\Big\} \\
&+ \beta \big{(}\sum_{i=1}^Nx_i\big{)}^2 +\sum_{i=1}^N h_i(x_i).
 \end{aligned}
\end{equation}
A simple arrangement using \eqref{g3} yields
\begin{equation}\label{g4}
\begin{aligned}
g(x)= - &\sum_{i=1}^N \min_{l_i \leq y_i\leq u_i} \Big \{ \beta
y^2_i + \Big{(}\beta  \sigma_{(-i)}(x)
 - \alpha  \Big{)}y_i + h_i(y_i)\Big\} \\
 & + \beta \big{(}\sum_{i=1}^N  x_i\big{)}^2 - \alpha^T x  + \sum_{i=1}^N h_i(x_i),
\end{aligned}
\end{equation}
where $\sigma^{(-i)}(x):= \sum_{j\neq i}^N x_j$.   From \eqref{g4}
 it follows  that evaluating $g(x)$, for each $x\in D$, one needs to
solve $N$-optimal problems each of them is  one-variable
minimization problem of the form
\begin{equation}\label{2.8}
\min_{l_i \leq y_i\leq u_i} \Big \{ \beta  y^2_i + \Big{(}\beta
\sigma^{(-i)}(x)
 - \alpha  \Big{)}y_i + h_i(y_i)\Big\}, \ i=1, 2, \ldots, N.
 \end{equation}
In order to compare the Cournot model presented above with existing models let us consider the Bertrand model. In a Bertrand model the firms producing a common homogenous commodity. In contrast to the Cournot model, here each firm sets prices rather than  the  production quantity. So, in  such a model, the demand is a function of price.and the customers buy from firms with lowest price. However, often this assumption is not realistic, since usually the products of the firms are not entirely interchangeable, and thus   some consumers may prefer one product to the other even it costs somewhat more.\\
\indent Suppose that the quantity level $x_i$   produced by firm $i$ depends  on the  price  and given by
$$x_i(p) = \gamma_i - \sigma_i p_i + \sum_{j\neq i}^n \lambda_{ij}p_j, \ i=1, \ldots, n$$
where $\gamma_i, \ \sigma_i > 0$,  $\lambda_{ij} \geq 0$ if $\ (j\neq i)$.  The condition   $\sigma_i >0$ means that  the  demand  for firm $i$ decreases as   its  price increases, while   $\lambda_{ij}  \geq 0$. ($i\not= j)$ means that the demand  for firm $i$ increases when other firms increase their price.\\
\indent The profit function of firm $i$ then is given as
$$f_i(p) := p_ix_i - h_i(x_i),$$
where, following \cite{Bi2}, we assume that  the cost  $h_i(.)$ is a concave function of the production level   and is given by
$$h_i(x_i) = \nu_ix_i - d_ix_i^2 \  \text{with} \ d_i \geq 0.$$
Then an elementary  computation   shows that the cost is a function of the price as
$$ \begin{array}{lll}
h_i(p) & =& -d_i\sigma_i^2p_i^2 + \sigma_i\big{[}2d_i \big{(} \gamma_i +  \sum_{j\neq i}^n\lambda_{ij}p_j\big{)}-\nu_i\big{]}p_i + \nu_i\big{(} \gamma_i +  \sum_{j\neq i}^n\lambda_{ij}p_j\big{)} \\
&& - d_i\big{(} \gamma_i +  \sum_{j\neq i}^n\lambda_{ij}p_j\big{)}^2
\end{array}$$
The profit function then  takes the form
$$\begin{array}{lll}
f_i(p) & = & \sigma_i(d_i\sigma_i - 1)p_i^2 + \big{[} \sigma_i \nu_i + (\gamma_i +  \sum_{j\neq i}^n\lambda_{ij}p_j)(1-2d_i \sigma_i)\big{]}p_i \\
& &+ d_i \big{(}\gamma_i +  \sum_{j\neq i}^n\lambda_{ij}p_j \big{)}^2 - \nu_i (\gamma_i +  \sum_{j\neq i}^n\lambda_{ij}p_j)
\end{array}$$

Each firm $i$ attempts to maximize its profit by choosing a corresponding price level on its strategy set $[0, T_i]$  by solving the optimization problem
$$f_i(p) = \max_{y_i \in [0,T_i]} f_i(p[y_i])), \ \ \forall i=1, \ldots, n,$$
where $p[y_i]$ is the vector obtained from $p$ by replacing $p_i$ with $y_i$.\\
\indent By the same technique as in the Nash-Cournot model the problem of finding a Nash equilibrium point of this Bertrand model can be formulated as a mixed variational inequality of the form
$$\begin{array}{lll} \text{Find} \  \ p\in T:= T_1\times\ldots\times T_n:
 \Phi(p,y):= \langle Gp - y, p-y \rangle + \psi(y)-\psi(p)\geq 0 \ \forall y\in T
\end{array}$$
 where
 $$ G = \left ( \begin{array} {ccccc}
0 & \lambda_{12}(1-2d_1 \sigma_1) & \ldots & \lambda_{13}(1-2d_1 \sigma_1)&\lambda_{1n}(1-2d_1 \sigma_1) \cr
\lambda_{21}(1-2d_2 \sigma_2) & 0 & \ldots & \lambda_{23}(1-2d_2 \sigma_2)&\lambda_{2n}(1-2d_2 \sigma_2) \cr
\ldots& \ldots& \ldots& \ldots&\ldots \cr
\lambda_{n1}(1-2d_n \sigma_n) & \ldots& \ldots&\lambda_{n,n-1}(1-2d_n \sigma_n)&  0
\end{array} \right )
$$
with
$$r_i = \gamma_i(1-2d_i \sigma_i), \ i=1, \ldots, n,$$
$$ \begin{array}{lll}
\psi(y) = \sum_{i=1}^n \sigma_i(d_i\sigma_i -1)y_i^2, \\
\end{array}
$$
So as the Nash-cournot model, the Bertrand model can be formulated as a mixed variational inequality $MV(D)$.
Note that since $\sigma_i(d_i \sigma_i -1), \ i=1, \ldots, n$  may be negative, the function  $\psi(.)$  may  not convex.
\section{An Algorithm for Global Equilibria}
In this section we describe an  algorithm  for  approximating a  global equilibrium point of the model. The idea of the proposed algorithm is quite natural, it uses
the convex envelope of the concave cost function to approximate the  original  model with the one having piecewise linear concave costs. The latter can be solved by an algorithm developed in \cite{MHQ} to obtain an approximate equilibrium point $x$. Then by evaluating  the gap function we can
check whether or not the obtained point $x$ is an $\epsilon$-equilibrium point. If not, we use an adaptive rectangular bisection   to get a better approximate point. Thanks to the rectangular structure of the strategy set and separability  of the cost function, the proposed algorithm can be implemented easily.
\subsection{A Search-Check-Branch Algorithm}
First we recall \cite{Ho1} that the convex envelope of a function $\varphi$ on a convex set $C$ is the convex function on $C$, denoted by $\co_C\varphi$ such that $\co_C\varphi(x) \leq \varphi(x)$ for every $x\in C$, and if $\xi$ is any convex function on $C$ satisfying $\xi(x) \leq \varphi(x)$ for
every $x\in C$, then $\xi(x) \leq \co_C\varphi(x)$ for every $x\in C$. It is well known \cite{Ho1} that the convex envelope of a concave function is affine, and that if $C= C_1\times \ldots\times C_N$ is  compact and $\varphi$ is separable, i.e., $\varphi(x_1,\ldots,x_N)= \sum_{j=1}^N \varphi_j(x_j)$ then $\co\varphi(x) = \sum_{j=1}^N\co\varphi_j(x_j)$ where $\co\varphi_j$ is the convex envelope of $\varphi_j$ over $C_j$. Clearly, since $h_i$, $i > n$ is affine,
$h_i \equiv coh_i$ on every convex set.\\
\indent The algorithm we are going to describe is a search-check-branch procedure. For a given tolerance $\epsilon \geq 0$, at each iteration,
  the algorithm consists of three steps. The  search-step requires solving convex quadratic programs for the approximate model with piecewise linear concave
 cost functions to obtain an approximate equilibrium point. The check-step  uses the gap function presented in the preceding section  to check whether the obtained solution is an $\epsilon$- equilibrium point or not yet. The branch-step employs  an adaptive rectangular bisection performed in the space of concave variables to obtain a better approximation for the model.\\
\indent To be precise, suppose that the strategy set $D :=  D_1\times \cdots\times D_N$. Let
$$I^0:=   D_1\times \ldots\times D_n, \ J^0:=  D_{n+1}\times \cdots\times D_N. $$
 For a $n$-dimensional subbox $I \subseteq I^0 $,  define
\begin{equation}\label{DI}
D_I:=\{x^T:= (x_1,\ldots, x_N): (x_1,\ldots,x_n) \in I, (x_{n+1},\ldots, x_N)\in  J^0 \}
\end{equation}
and consider the convex mixed variational inequality CMV($D_I$) defined as
$$ \begin{array}{l}
\text{Find} \ x^{D_I}\in D_I \ \text{such that:}\\
 \langle \tilde{B}_1x^{D_I}- \alpha, y-x^{D_I} \rangle + y^TB_1y + \co_{I}h(y)
  + \sum_{j=n+1}^N h_j(y)\\ - (x^TB_1x^{D_I} + \co_{I}h(x^{D_I})+ \sum_{j=n+1}^N h_j(x^{D_I})
) \geq 0\ \forall y\in D_I.
\end{array}
$$
In what follows we write $x^{D_I}= ( x^I,x^J)$ with $x^I \in I$,
$x^J\in J^0$.\\
\indent Since $\co_{I}h(.)$ is affine, by Lemma \ref{2.1a},  this problem is reduced to the strongly convex quadratic program
$$\min_{x\in D_I}\{ x^TQx + (c^I)^Tx\}, \eqno(QD_I)$$
where $Q:= \dfrac{1}{2}\tilde{B_1} + B_1$, $c^I = (c^{I_1},  \ldots,
c^{I_N})^T$ with $c^{I_j}:= (a^{I_j} - \alpha)     (j=1, 2, \ldots,N)$.\\
\indent Suppose that each strategy set $D_j$ ( $j=1,\ldots,n$)  has been divided into interval $D_{j,1},\ldots, D_{j,k_j}$  on each of
  them the cost function is affine. Let  $\Delta$ be the  set of $n$-dimensional subboxes defined as
  $$\Delta:=\{ B:= I_1\times\cdots\times I_n: I_j \in \{D_{j,1}, \ldots , D_{j, k_j}\}, j = 1,\ldots,n\}.$$
Define $\Sigma$ as the family of $N$-dimensional subboxes by taking
$$\Sigma:=\{ I= B\times J^0: B \in \Delta\}.$$
\indent Let us define the gap function for the model with piecewise concave cost function, that is
\begin{equation}\label{gapb}
 \bar{g}(x):= - \min_{y\in D}\bar{\phi}(x, y)
\end{equation}
where
\begin{equation}\label{phib}
\bar{\phi}(x, y)  := \langle \tilde{B_1}x - a, y - x \rangle
 + y^TB_1y - x^TB_1x + \bar{h}(y) -
 \bar{ h} (x)\},
\end{equation}
where $\bar{h}$ is the piecewise linear concave function obtained by
taking the convex envelope of $h$ on each  element  of $\Sigma$.\\
\indent Note that, since $h_i$ is affine on $D_i$ for every $i=n+1,\ldots,N$, the convex envelope of $h_i$ on any subbox coincides
with $h_i$ for every $i= n+1,\ldots,N$. In particular, $\co_D h$ is affine and
$$\co_D h(x) = \sum_{j=1}^n\co_{I^0} h_j(x) + \sum_{i=n+1}^N h_i(x).$$
\indent First we briefly describe the algorithm in \cite{MHQ} as follows.\\
\noindent {\bf Algorithm 1} (Search-and-Check). Choose a tolerance $\epsilon \geq 0$.\\
\noindent {\bf Step 1}: Select a subbox $I\in  \Sigma$.\\
\noindent {\bf Step 2}: Solve the strongly convex quadratic problem $(QD_I$) to obtain its unique solution $x^{D_I}$.\\
\noindent {\bf Step 3}:\\
\indent a) If $\bar{g}(x^{D_I}) \leq \epsilon$, terminate: $x^{D_I}$ is an $\epsilon$-equilibrium point for piecewise concave cost model.\\
\indent (b) If $\bar{g}(x^{D_I}) > \epsilon$ and $\Sigma = \emptyset$, then terminate: the model has no equilibrium point. Otherwise, replace
$\Sigma$  by $\Sigma \setminus \{I\}$ and return to {\bf Step 1}.\\
\indent It is obvious that  in the worst case, the algorithm searches all subboxes in $\Sigma$, however the computational results reported in \cite{MHQ} show that by  using the gap function, in general, the algorithm finds an $\epsilon$- equilibrium point without searching all elements of $\Sigma$.\\
\indent Using  Algorithm $1$ described above we can develop an algorithm for approximating an equilibrium point of the model where some of  the cost functions are concave. The idea is quite natural. In fact, at each iteration we use the convex envelope of the concave cost function to obtain a model with
piecewise lineae concave cost function to which we can apply the search-and-check Algorithm $1$ to obtain an approximate equilibrium point. If the obtained point is not yet an $\epsilon$- equilibrium point, we use an adaptive rectangular bisection (Rule $1$ below)  to reduce the difference between the concave function and its convex envelope to obtain  a better approximate equilibrium point for the original model, and so on.\\
\noindent {\bf An adaptive rectangular  bisection} (Rule 1). Let $I$ be a given $n$-dimensional  subbox of $D_1\times\ldots\times D_n$. For $x^I\in I$, define
$$j_{max}:= \argmax_{1\leq j \leq n} \{ h_j(x^I_j)-\co h_j(x^I_j)\}.$$ Then we bisect $I$ into two boxes  via the middle point  of edge $I_{j_{max}}$. We call this middle point  the {\it bisection point} and $ j_{\max}$ the {\it bisection index}.\\
\indent For this bisection we have the following lemma whose proof can be found, e.g., in \cite {MO1,Mu1}.
\begin{lemma}\label{bi} Let $\{I^k\}$  be an infinite sequence of  boxes generated by the adaptive rectangular bisection  Rule $1$ such  that $I^{k+1}\subset
I^k$ for every $k$. Let  $b^k$ be the bisection point and $j_k$  be the bisection index for $I^k$. Then $lim_{k\to \infty}( h_{j_k}(b^k)- \co_{I^k}h_{j_k}(b^k) ) = 0$.  Consequently, $\{I_{j_k}\}$ tends to a singleton. provided $h_{j_k}$ is (concave) not affine on $I_{j_k}$ for every $j_k$.
\end{lemma}
\indent For each subbox $I$ having $n$-edges $I_j$ $(j=1,\ldots,n)$ we define
$$\rho(I_j):= \max_{t\in I_j} \{ h_j(t) - coh_j(t)\}$$
and
\begin{equation}\label{roI}
\rho(I):=\max\{ \rho(I_j): \j=1,\ldots,n\}.
\end{equation}
The algorithm now can be described as follows:\\
\noindent {\bf Algorithm 2} (Search-Check-Branch for global equilibria).\\
\noindent {\it Initial step}. Choose a tolerance $\epsilon \geq 0$, take the initial box $I^0 := D_1\times\ldots\times D_n$. Solve the convex mixed
variational inequality CMV($D$) defined as
$$\text {Find}\ x \in D: \overline{\Phi}_0(x,y):= \langle \tilde{B_1}x - a, y - x \rangle
 + y^TB_1y - x^TB_1x +   \co_D h(y) -
 \co_D h(x)\geq 0\forall y\in D,$$
 which is equivalent to the strongly convex quadratic program $(QD_{I^0})$ to obtain its unique  solution $u^0$.\\
\indent Let $\Sigma_0:= \{I^0\}$. $x^0 := u^0$.\\
\noindent {\bf Iteration $k$ } ($k= 0,1\ldots$)\\
\indent At the beginning of each iteration $k$ we have:\\
\indent $\bullet$ $\Sigma_k$: a finite family of $n$-dimensional subboxes of $ I^0$;\\
\indent $\bullet$ $ u^k = (u^{k_1},u^{k_2})$ with $u^{k_1} \in I^0$, $u^{k_2} \in J^0$, the equilibrium point of the model with piecewise linear
 concave function;\\
\indent $\bullet$  $x^k\in D$: the currently best feasible point, i.e., $g(x^k)$ is smallest among the obtained
 feasible points so far.\\
\noindent {\bf Step 1}.\\
\indent a) If $g(x^k) \leq \epsilon$, terminate: $x^k$ is an $\epsilon$-equilibrium point of the original model.\\
\indent b) If $g(x^k)  > \epsilon$, choose $I^k \in \ \Sigma_k$ such that
$$ \rho(I^k) = \max \{ \rho(I); \ I\in \Sigma_k\}. $$
\noindent {\bf Step 2}. Use the bisection Rule 1 described above to bisect $I^k$ into two
boxes $I^{k^+}$ and $I^{k^-}$. Let $j_k$ be the bisection index for $I^k$.\\
\noindent {\bf Step 3}. Solve the strongly convex quadratic program $(QD_I)$
with $I = I^{k^-}$ and $I= I^{k^+}$ to obtain $x^{k+}$ and $x^{k-}$ respectively.\\
\noindent {\bf Step 4}. If either $g(x^{k+}) \leq \epsilon$ or  $g(x^{k-})\leq \epsilon$, terminate.\\
\noindent Otherwise, update $x^k$,  $\Sigma_k$ and the  linear piecewise concave cost function by taking respectively
$$x^{k+1} \in \{x^k, x^{k+}, x^{k-}\} \ \text{ such that}\ g(x^{k+1}) = \min \{g(x^k), g(x^{k-}), g(x^{k+})\}, $$
$$\Sigma_{k+1} =(\Sigma_k \setminus  \{ I^k\}) \cup\{I^{k^-}, I^{k^+}\}.$$
\noindent {\bf Step 5}. Compute the convex envelope of function $h_{j_k}$ on the egde $j_k$ of the subboxes $I^{k^-}$, $I^{k^+}$, thereby to obtain the new approximation bifunction
$$\overline{\Phi}_{k+1}(x,y):= \langle \tilde{B_1}x - a, y - x \rangle
 + y^TB_1y - x^TB_1x +   \co_{k+1} h(y) -
 \co_{k+1} h(x),$$
 where $\co_{k+1}h$ is the convex envelope of $h$ obtained by replacing the convex envelope of $h_{j_k}$ on the edge $j_k$ of $I^k$ by the convex envelope of $h_{j_k}$ on the edge $j_k$ of $I^{k^-}$ and $I^{k^+}$. Then use Algorithm 1 with the just obtained piecewise linear concave cost function to solve the newly approximated piecewise  linear concave model to obtain $u^{k+1}$.\\
Increase $k$ by one and  go to {\bf Step 1} of iteration $k$.\\
\noindent Suppose that  every model with piecewise linear concave cost function has an $\epsilon$- equilibrium point for any $\epsilon > 0$.  Then we have the
following convergence result.\\
\noindent {\bf  Convergence Theorem.}\\
{\it \indent (i)  If the algorithm terminates at iteration $k$ then  $x^k$ is an $\epsilon$-equilibrium point.\\
\indent (ii) If the algorithm does not terminate, it generates an infinite sequence $\{ x^k\}$ such that any its cluster point is an equilibrium point whenever the model has an equilibrium point.
Furthermore $g(x^k) \searrow 0$ as $k\to \infty$.}\\
\noindent {\it Proof.} The statement (i) is obvious.\\
\noindent To prove statement (ii) we suppose that the algorithm never terminates. Let $x^*$ be any cluster point of $\{x^k\}$. Then
there exists a subsequence  of $\{x^{k_q}\}$ that tends to $x^*$. Thus the corresponding sequence of selected intervals has a nested sequence, which,  by taking a subsequence if necessary, we denote also by $I^{k_q}$. Since $I^{k_q}$ is the box to be bisected at iteration $k_q$, by Lemma \ref{bi}, $\{I^{k_q}\}$ tends to a singleton, which implies that $h_{j_q}(x_{j_q}) - \co h_{j_q}(x_{j_q}) \to 0$ as $q\to \infty$  ($j_q$ is the bisection index at iteration $k_q$). By the rule for selecting the bisection index, we have $h_j(x_j ) - \co h_j(x_j ) \to 0$ for every $j$. Since $u^{k_q}$ is an  equilibrium point of the model with piecewise linear concave cost function, we have $\bar{g}_{k_q}(u^{k_q})    = 0$ for every $q$, where $\bar{g}_{k_q}$ is the gap function for the piecewise linear concave cost model at iteration $k_q$. By the definition of the gap function $g$ for the original model and of $\bar{g}$ for the approximate model, and the rule for selecting bisection index, we can write
$$\bar{g}(u^{k_q}) - 2\sigma_{k_q} \leq g(u^{k_q}) \leq \bar{g}(u^{k_q}) + 2\sigma_{k_q}   \  \forall q.$$
Letting $q\to \infty$, since $\sigma_{k_q} \to 0$, $u^{k_q} \to u^*$, by continuity of $g$, we obtain $g(u^*) = 0$.\\
On the other hand, since $x^{k_q}$ is the currently best feasible point obtained at iteration $k_q$, we have $0 \leq g(x^{k_q}) \leq g(u^{k_q})$. Letting $q \to \infty$, by continuity of $g$, we obtain $0 \leq g(x^*) = g(u^*) = 0$, which means that $x^*$ is an equilibrium of the model.  Note that, since $x^k$ is the currently best feasible point obtained at iteration $k$, by definition, the sequence $\{g(x^k)\}$ is nonincreasing. Since the whole sequence $\{x^k\}$ is bounded, it has a subsequence $\{x^{k_j}\}$ converging to some $\bar{x}$. Then, as we just have shown, $\bar{x}$ is an equilibrium point which implies $g(\bar{x}) = 0$. Then the whole sequence $\{g(x^k)\}$ tends to $0$ as well. $\hfill \square$
\begin{Remark}\label{4.2}
In order to save the memory, we may use a criterion to delete every subbox that does not contain  an equilibrium point in it.
\end{Remark}
 The following lemma gives a criterion that can be used to check whether  a  subbox
 contains an equilibrium  point or not.  In fact,  for a  subbox $D_I:= \{x\in D: l^I \leq x\leq u^I\}$,
  let us define the numbe $$\tilde{g}(D_I):= - \min_{y\in D_I}\{\langle \tilde{B}_1u^I - a, y
 \rangle +y^TB_1y +h(y)\} - (l^I)^T B_1 l^I + a^T u^I - h(l^I). $$
Then we have the following lemma:
 \begin{lemma} \label {del} Suppose $x^{D_I}$ is an optimal solution of Problem $(QD_I)$.\\
\indent (i) If $\co_Ih(x^{D_I}) = h(x^{D_I})$ then $x^I$ is the  equilibrium point the model restricted on $D_I$.\\
\indent (ii)  If $\tilde{g}(D_I) > 0$,  the subbox $D_I$ contains no
equilibrium point of the model.
\end{lemma}
{\it Proof.}\\
\indent (i) Since $x^{D_I}$ is the solution of $(QD_I)$,  we have
$$\langle \tilde{B}_1x^{D_I} - a, y-x^{D_I} \rangle + y^TB_1y + \co_Ih(y) - (x^I)^TB_1x^{D_I} - \co_Ih(x^{D_I}) \geq 0, \forall y\in D_I.$$
  Note that $h(y) \geq \co_Ih(y), \ \forall y\in D_I$, by the assumption, $\co_Ih(x^{D_I}) = h(x^{D_I})$, we obtain
$$\langle \tilde{B}_1x^{D_I} - a , y-x^{D_I} \rangle + y^TB_1y + h(y) -
 (x^I)^TB_1x^{D_I} - h(x^{D_I}) \geq 0$$
for every $y\in D_I$,  which means that $x^{D_I}$ is the equilibrium point the model restricted on $D_I$.\\
\indent (ii) We now prove that $g(x)  > 0 $ for all $x\in D_I$. Indeed,  by definition
$$\Phi(x, y) = \langle \tilde{B}_1x - a, y \rangle +y^TB_1y +h(y) - x^TB_1x + a^T x - h(x).$$
 Since $y\geq 0, \tilde{B}_1$ and $B_1$ are non-negative matrices, $h_i(.) (i=1, 2, \ldots, n)$ are increasing functions and $l^I \leq x \leq u^I $  for every $x \in D_I$,   we can  write, for every $y\in D$ and $x\in D_I$,
\begin{equation}\label{gg}
\begin{aligned}
\Phi(x,y) & =  \langle \tilde{B}_1x - a, y \rangle +y^TB_1y +
h(y) - x^TB_1x + a^T x - h(x)\\
& \leq  \langle \tilde{B}_1u^I - a, y \rangle +y^TB_1y +h(y) -
(l^I)^TB_1l^I + a^Tu^I - h(l^I) .
\end{aligned}
\end{equation}
By the definition of  $\tilde{g}(D_I)$, it follows from \eqref{gg} that
$$g(x): = - \min_{y\in D} \phi(x,y))\geq - \min_{y\in D_I} \phi(x,y))  \geq \tilde{g}(D_I)> 0 \  \forall\   x\in D_I,$$
which implies that $D_I$ does not contain an equilibrium point.  $\hfill \square$
\section{An Algorithm for Local Equilibria }
Using  the fact that a point $x^*\in D$ is an equilibrium point of the model if and only if the gap function $g(x^*) = 0$, we
 say that a point $\bar{x}$ is a {\it local equilibrium point} of the model if there exists an open set $B \subset D$ such that $\bar{x}\in B$, $g_B(\bar{x}) = 0,$ where  $g_B$ stands for the gap function of the model restricted on $B$. Note that because of concavity of the cost function, in this equilibrium Nash-Cournot model, a local equilibrium point may not be a global one.\\
\indent In this section, we propose an algorithm for approximating a local equilibrium point of the model by using again the  gap function.\\
Namely,  for a subbox $$I:=\{ x=(x_1,\ldots,x_n)^T : l_i \leq x_i \leq
u_i, \ i=1,\ldots,n\}, $$  let, as before,  $D_I$ be the subbox of $D$
 consists of all points  $x^T =(x_1,\ldots, x_n,\ldots, x_N)$ such that $(x_1,\ldots,x_n) \in I$. That
 is
$$D_I =\{ x=(x_1,\ldots,x_N)^T, \ l_i\leq x_i \leq u_i, i=1,\ldots N\}.$$
Then define the gap function $g_{D_I}$ restricted on $D_I$ by taking
\begin{equation}\label{glo}
\begin{aligned}
g_{D_I}(x)= - &\sum_{i=1}^N \min_{l_i \leq y_i\leq u_i} \Big \{
\beta y^2_i + \Big{(}\beta \sigma_{(-i)}(x)
 - \alpha \Big{)}y_i + h_i(y_i)\Big\} \\
 & + \beta \sum_{i=1}^N  x^2_i - a^T x  + \sum_{i=1}^N h_i(x_i),
\end{aligned}
\end{equation}
where $\sigma^{(-i)}(x):= \sum_{j\neq i}^N x_j$.
As before we use the convex envelope of the concave function $h$ on each  subbox $D_I$ to obtain a convex mixed variational inequality
whose solution can be obtained by solving a strongly convex quadratic over $D_I$. If it happens that at the obtained solution
the values of the  cost function   and its convex envelope on $D_I$ coincide, this solution is a local equilibrium point of the model.
Otherwise we bisect $I$ to reduce the difference between  the cost function and its convex envelope on $ D_I$. Note that if
$g_{D_I}(x) = 0$ for some $x\in D_I$, then $x$ is a local equilibrium point. Thus, if $x\in D_I$ and $g_I(x) \leq \epsilon$, then $x$ is an $\epsilon$-local equilibrium point. Since $x^{D_I}$ is the equilibrium point of the model with respect
to $D_I$, from the definitions of the convex envelope of $h$  and the gap function restricted on $D_I$, it follows that $h(x^{D_I}) -
\co_I h(x^{D_I}) = 0$ implies $g_{D_I}(x^{D_I}) = 0$.  The algorithm now can be described as follows.\\
\noindent{\bf Algorithm 3} (Search-Check-Branch for local equilibria).
\noindent{\it Initial step}. Choose  tolerances $\epsilon > 0$  and solve Problem  (Q$D$)
 to obtain its optimal solution  $x^{I^0}$.\\
\noindent Compute $\rho_0 := \rho(I^0)$ and  $\epsilon_0:=  g_D(x^{I^0})$. Set the initial box $I^0$ and let
$\Gamma_0 :=\{I^0\}$.\\
\noindent {\it Iteration} $k \ (k=0, 1, \ldots).$ At the beginning of each iteration $k$ we have:\\
\indent $\bullet$ $\Gamma_k$:  a finite family  of $n$-dimensional subboxes of $I^0$;\\
\indent $\bullet$ $\epsilon_k = \min \{ g_{D_I}(x^{D_I}): I\in \Gamma_k\}$, where $x^{D_I}$ is the optimal solution of the convex quadratic
program (Q$D_I)$;\\



\noindent {\bf Step 1.} (Stoping criteria) If  $\epsilon_k \leq  \epsilon$, terminate: $x^{D_I}$ with $\epsilon_k = g_{D_I}(x^{D_I})$ is an $\epsilon$-local equilibrium point.


\noindent {\bf Step 2.} (Selection) Choose $I^k \in \Gamma_k$ such that
$$\rho_k := \rho(I^k) = \max\{\rho(I):  I\in \Gamma_k\}.$$
 \noindent {\bf Step 3.} (Bisection): Divide the subbox $I^k$ into
two subboxes $I^{k^+}$ and $I^{k^-}$ by the  bisection Rule $1$.\\
\noindent {\bf Step 4.} Solve the strongly convex quadratic programs (Q$D_I$)
with $I:= I^{k^+}$
  and  $I:= I^{k^-}$ to obtain the optimal solutions
  $x^{k^+}$ and $x^{k^-}$ respectively.
 Compute  $\rho(I^{k^+})$ and $\rho(I^{k^-})$.\\
\noindent {\bf Step 5.} Let $ \epsilon_{k+1}:= \text{argmin} \{\epsilon_k,
g_{D_I}(x^{D_I}) \ \text{with}\  I = I^{k^-}\  \text{and} \ I = I^{k^+} \}.$\\
\noindent {\bf Step 6.} (Updating) If $g(x^{D_I}) > 0$ delete $I$ from
 further consideration.\\
\indent Let $\Gamma_{k+1}$ be the remaining
 set. If $\Gamma_{k+1} = \emptyset$, terminate: the model has no
 $\epsilon$-local equilibrium point. Otherwise,  go to iteration $k$ with
 $k:=k+1$.\\
\noindent {\bf  Convergence.} {\it The algorithm terminates after a finite iteration yielding an $\epsilon$- local equilibrium point whenever
it does exist}.\\
\indent The proof of this convergence result is evident because of the fact that $\epsilon > 0$,  that the sequence of selected boxes tends to a
singleton and that the gap function is continuous.
\begin{Remark}\label{4.3} If for every $i$, the cost function $h_i$ satisfies the condition
 \begin{equation}\label{ex}
 \nabla^2h_i(y_i) \geq - 2\beta_i \ \forall y_i \in D_i.
 \end{equation}
Then the model admits a solution.\\
Indeed, for each $x\in D$,  let   $H_i(x)$ be the solution set of the problem
$$\min_{y_i\in D_i}\Big {\{} \varphi_i(x_{-i},y_i):= \beta y_i^2 + \Big{(} \beta \sum_{j\neq i}^Nx_j - \alpha_i\Big{)}y_i +h_i(y_i)\Big{\}}. \eqno(QD_i(x)).$$
 It is easy to check that condition \eqref{ex} ensures that the object function of this problem is convex in $y_i$. Thus $H_i(x)$ is a closed convex of the interval $D_i$. Since the objective function of this problem is continuous and the feasible is compact, the solution set $H_i(x)$ is a upper semicontinuous mapping from $D$ into itself, by well-known Kakutani fixed point, the mapping $H(x):= H_1(x)\times \ldots\times H_N(x) $ has a fixed point $x^*$, which is also an equilibrium point of the model.\\
\indent Note that both the cost functions
$$h_i(y_i) = \ell_iy_i - d_i y_i^2$$
with $\beta_i > d_i > 0$ for all $i=1, \ldots, n$.
used in \cite{Bi2} and
$$h_i(y_i) = \mu_iy_i + \ln(1+ \gamma_iy_i),$$
 with $\gamma_i > 0$   and $\gamma_i^2 \leq 2\beta_i, \ \forall i = 1, \ldots, n$
satisfy condition \eqref{ex}.
\end{Remark}
\section{Computational Results and Experiments} The   proposed  two  algorithms
    were implemented in MATLAB. The programs were executed on a
PC Core 2Duo 2*2.0 GHz, RAM 2GB. We   tested the program  on
different groups of problems, each of them contains ten problems of
 different sizes $N$ and $n$, but having  randomly generated input data.
Namely, for each problem, the numbers  $\alpha$, $\beta$, $\mu_i$
($i=n+1,\ldots,N)$ are randomly generated in the interval $[ 20. 30] $,
$[0.001, 0.005]$ and  $[10. 20]$ respectively. We take the  cost
functions of the form
\begin{equation}\label{excost}
h_j(x_j) = a_j x_j + ln(1+ \gamma_j x_j) , \ (j=1,\ldots,n), \
h_i(x_j):= \mu_i x_i \ (i = n+1,\ldots,N).
\end{equation}
 where $\gamma_j  $ and   $a_j$   are  randomly generated in
 $[7,15]$  and  $[2,7]$    respectively. The strategy set of firm $i$ is
 $D_i := [0, u_i]$ where
 each  $u_i$ is randomly generated
in the interval $[100. 500]$.

 The obtained results are reported  in Table $4.1$ below,  where we use the following headings:
\begin{itemize}
\item $N$:    number of the firms;
   \item $n$:    number of the  firms  having  concave (but not affine) cost;
       \item {\it Average time}:    the average  time (in second) needed to solve one
     problem;
   \item  {\it Average iter}:  the average numbers of
     iterations for   one problem.
\item  {\it Glob-GSCB}: number of problems for which an equilibrium point
 was obtained by  Search-Check-Branch Algorithm for global equilibria.
  \item {\it Glob-LCB}:  number of problems for which  a global
       optimal solution was obtained by  Search-Check-Branch for local equilibria.
\end{itemize}

\begin{center}
\begin{tabular}{|c|c|p{1.5cm}|p{1.5cm}|p{1.5cm}|p{1.5cm}|p{1.5cm}|p{1.5cm}|}
\hline
 \multicolumn{2}{|c|}{Size} & \multicolumn{3}{|c|}{GSCB-Alg.} &
\multicolumn{3}{|c|}{LSCB-Alg.} \\ \hline N & n & Average time
& Average iter. &Glob-GSCB & Average time & Average iter & Glob-LSCB \\
\hline
5 & 5 & 0.00 & 1 & 10 & 0.03 &1& 10 \\
50 & 5 & 8.98 & 133 & 10 & 0.06 & 1&10 \\
100 & 5 & 17.89 & 171 & 10 & 0.18 & 2&8 \\
200 & 5 &1.78 & 7 & 10 & 0.29 & 2&8 \\  \hline
10 & 10 & 9.65 & 308 & 10 & 0.05 &1& 8 \\
50 & 10 & 82.35 & 1141 & 10 & 0.22 & 4 & 4 \\
100 & 10 & 47.05 & 445 & 10 & 0.43 & 5& 7 \\
200 & 10 & 41.06 & 203 & 10 & 0.33 & 2& 7 \\ \hline
20 & 20 & 127.15 & 2478 & 10 & 1.29 & 24&1\\
50 & 20 & 98.10 & 1231 & 10 & 0.50 & 7 &3\\
100 & 20 & 105.00 & 914 & 10 & 1.72 & 16 &3 \\
200 & 20 & 440.88 & 2216 & 10 & 1.94 & 11& 5 \\ \hline
30 & 30 & 286.57 & 3754 & 10 & 0.89 & 13 & 2 \\
50 & 30 & 246.44 & 2901 & 10 & 1.23 & 17& 1\\
100 & 30 & 872.27 & 7193 & 10 & 0.73 & 7 & 2 \\
200 & 30 & 750.72 & 3514 & 10 & 2.70 & 15 & 4 \\ \hline
40 & 40 & 515.10 & 5944 & 10 & 3.09 & 40 & 2 \\
50 & 40 & 1332.10 & 14820 & 9 & 7.69 & 97 & 0 \\
100 & 40 & 646.53 & 5213 & 10& 2.85 & 26 & 0 \\
200 & 40 & 898.09 & 4169 & 9& 3.83 & 21 & 1  \\ \hline
100 & 100 & Skip &  - &  - & 20.21 & 148 & 0 \\
200 & 100 & Skip &  - &  - & 132.64 & 568 & 0 \\
200 & 200 & Skip &  - &  -& 107.63 & 400 & 0 \\
300 & 200 & Skip &  - &  - & 252.67 & 579 & 0 \\ \hline
\end{tabular}
\end{center}
\centerline{Table 4.1}
\indent From the obtained results reported in Table $4.1$ we can conclude the
followings for the tested concave cost functions given as \eqref{excost}.\\
\indent $\bullet$ Algorithm $2$ for global equilibrium point can solve
models with a moderate number ($n\leq 40)$ of concave cost
functions, while Algorithm $3$ can solve models where the number of
concave cost functions much larger.\\
\indent $\bullet$  For models where the   number of the firms having concave
cost is somewhat   large ($n \geq 40$), the   local equilibrium
point obtained by the local algorithm is often not a global one.
\section{Conclusion}
A Nash-Cournot oligopolistic equilibrium
model involving  concave cost functions may have local equilibrium
points that are not global ones. We have approximated such a model
with the one having piecewise linear concave function by using the
convex envelope of a separable concave function over a box. Based
upon this approximation we have proposed two algorithms for
approximating a global as well as local equilibrium points that
employ a gap function as a stoping criterion for the algorithms, and
an update rectangular bisection to make the approximation better.
Some computational results have been reported showing efficiency of
the proposed algorithms for models where the number of the concave
(but not affine) cost functions is not  large ($n\leq 40$) for
global algorithm, and $(n \leq 200$) for local one. An open question
that  would be interesting for further consideration is to find a
differentiable gap function, for which a local optimization
algorithms such as descent    ones in \cite{Fu1a} or DCA in
\cite{Ph1} could be applied efficiently.
\section*{Acknowledgements}
This work is supported by the National Foundation for Science and Technology Development (NAFOSTED), Vietnam.

\end{document}